\documentclass[12pt]{article}
\usepackage{mathrsfs}

\usepackage{amsmath, epsfig, cite, lineno}
\usepackage{amssymb,amsthm}
\usepackage{amsfonts, color}
\usepackage{latexsym}

\newtheorem{thm}{Theorem}[section]

\newtheorem{cor}[thm]{Corollary}
\newtheorem{conj}[thm]{Conjecture}

\theoremstyle{remark}
\newtheorem{remark}[thm]{Remark}

\def\N{{\mathbb N}}

\DeclareMathOperator{\podd}{pod}

\numberwithin{equation}{section}

\begin{document}


\begin{center}
{\Large\bf Two truncated identities of Gauss}
\end{center}

\vskip 2mm \centerline{Victor J. W. Guo$^1$  and Jiang Zeng$^{2}$}
\begin{center}
{\footnotesize $^1$Department of Mathematics, East China Normal University,\\ Shanghai 200062,
 People's Republic of China\\
{\tt jwguo@math.ecnu.edu.cn,\quad http://math.ecnu.edu.cn/\textasciitilde{jwguo}}\\[10pt]
$^2$Universit\'e de Lyon; Universit\'e Lyon 1; Institut Camille
Jordan, UMR 5208 du CNRS;\\ 43, boulevard du 11 novembre 1918,
F-69622 Villeurbanne Cedex, France\\
{\tt zeng@math.univ-lyon1.fr,\quad
http://math.univ-lyon1.fr/\textasciitilde{zeng}} }
\end{center}


\vskip 0.7cm \noindent{\bf Abstract.}
Two  new expansions for partial sums of Gauss' triangular and
square numbers series are given. As a consequence, we derive a family of inequalities for
 the overpartition function $\overline{p}(n)$ and for the  partition function $\podd (n)$ counting the partitions of $n$ with distinct odd parts.
Some further  inequalities for variations of partition function are proposed as conjectures.



\section{Introduction}
The partition function $p(n)$ has the generating function
\begin{align*}
\sum_{n=0}^{\infty} p(n) q^n
=\prod_{n=1}^\infty \frac{1}{1-q^n}=1+q+2q^2+3q^3+5q^4+7q^5+11q^6+\cdots.  
\end{align*}
Two classical results in the partition theory \cite[p.~11]{Andrews98}  are  Euler's pentagonal number theorem
\begin{align}
1+\sum_{j=1}^{\infty} (-1)^j (q^{j(3j-1)/2}+q^{j(3j+1)/2})
=\prod_{n=1}^\infty (1-q^n), \label{eq:euler2}
\end{align}
and Euler's  recursive  formula for computing
$p(n)$:
\begin{align}\label{eq:recur1}
p(n)+\sum_{j=1}^{\infty} (-1)^{j}(p(n-j(3j-1)/2)+p(n-j(3j+1)/2))= 0,
\end{align}
where $p(m)=0$ for all negative $m$.

Recently, Merca \cite{Merca} stumbled upon the following inequality:
\begin{align}
p(n)-p(n-1)-p(n-2)+p(n-5)\leqslant 0, \label{eq:merca}
\end{align}
and then,
Andrews and Merca \cite{AM} proved more generally  that, for $k\geqslant 1$,
\begin{align}
(-1)^{k-1}\sum_{j=0}^{k-1} (-1)^{j}(p(n-j(3j+1)/2)-p(n-(j+1)(3j+2)/2))\geqslant 0 \label{eq:ineq1}
\end{align}
with strict inequality if $n\geqslant k(3k+1)/2$.

The $q$-shifted factorial and $q$-binomial coefficient  are defined by
\begin{align*}
(a)_\infty=(a;q)_\infty=\prod_{n=0}^{\infty}(1-aq^n),
\quad (a)_M=(a;q)_M=\frac{(a;q)_\infty}{(aq^M;q)_\infty},
\end{align*}
and
\begin{align*}
{M\brack N}={M\brack N}_q=\frac{(q;q)_M}{(q;q)_N (q;q)_{M-N}}.
\end{align*}
Whenever the base of a $q$-shifted factorial or $q$-binomial
coefficient is just $q$ it will be omitted.
The proof of \eqref{eq:ineq1} in \cite{AM} is based on the truncated formula of \eqref{eq:euler2}:
\begin{align*}
\frac{1}{(q)_\infty}\sum_{j=0}^{k-1} (-1)^j q^{j(3j+1)/2}(1-q^{2j+1})
=1+(-1)^{k-1} \sum_{n=k}^\infty
\frac{q^{(k+1)n+{k\choose 2}}}{(q)_n}{n-1\brack k-1}.
\end{align*}

Motivated by Andrews and Merca's work \cite{AM}, in this paper we shall prove
new truncated forms of two identities of Gauss \cite[p.~23]{Andrews98}:
\begin{align*}
1+2\sum_{j=1}^{\infty} (-1)^j q^{j^2}
&=\frac{(q)_\infty}{(-q)_\infty},  \\  
\sum_{j=0}^\infty (-1)^jq^{j(2j+1)}(1-q^{2j+1})
&=\frac{(q^2;q^2)_\infty}{(-q;q^2)_\infty}, 
\end{align*}
and  derive similar overpartition function and special partition function inequalities.
\begin{thm}\label{thm:1}
For $|q|<1$ and $k\geqslant 1$, there holds
\begin{align}
\frac{(-q)_\infty}{(q)_\infty}\bigg(1+2\sum_{j=1}^k (-1)^j q^{j^2}\bigg)
=1+(-1)^k \sum_{n=k+1}^\infty \frac{(-q)_{k} (-1)_{n-k}q^{(k+1)n}}{(q)_n}{n-1\brack k}. \label{eq:genovp}
\end{align}
\end{thm}
The overpartition function $\overline{p}(n)$, for $n\geqslant 1$, denotes the number of ways of writing the integer $n$ as a sum of
positive integers in non-increasing order in which the first occurrence of an integer may be
overlined or not, and $\overline{p}(0)=1$
(see Corteel and Lovejoy\cite{CL}).  It is easy to see that
\begin{align}
\sum_{n=0}^{\infty}\overline{p}(n) q^n
=\frac{(-q)_\infty}{(q)_\infty}
=1+2q+4q^2+8q^3+14q^4+24q^5+40q^6+\cdots. \label{eq:ovp}
\end{align}

\begin{cor}\label{cor:2} For $n,k\geqslant 1$, there holds
\begin{align}
(-1)^k \bigg(\overline{p}(n)+2\sum_{j=1}^k (-1)^{j}\overline{p}(n-j^2) \bigg)\geqslant 0  \label{eq:neq-overp}
\end{align}
with strict inequality if $n\geqslant (k+1)^2$. For example,
\begin{align}
\overline {p}(n)-2\overline {p}(n-1)\leqslant 0,  \label{eq:neq-overp2} \\
\overline {p}(n)-2\overline {p}(n-1)+2\overline {p}(n-4)\geqslant 0, \nonumber \\
\overline {p}(n)-2\overline {p}(n-1)+2\overline {p}(n-4)-2\overline {p}(n-9)\leqslant 0. \nonumber
\end{align}
\end{cor}

\begin{thm}\label{thm:3}
For $|q|<1$ and $k\geqslant 1$, there holds
\begin{align*}
&\hskip -2mm \frac{(-q;q^2)_\infty}{(q^2;q^2)_\infty}
\sum_{j=0}^{k-1}(-1)^jq^{j(2j+1)}(1-q^{2j+1}) \nonumber\\
&=1+(-1)^{k-1} \sum_{n=k}^\infty \frac{(-q;q^2)_{k} (-q;q^2)_{n-k}q^{2(k+1)n-k}}
{(q^2;q^2)_n}{n-1\brack k-1}_{q^2}.
\end{align*}
\end{thm}
Following Hirschhorn and Sellers \cite{HS} we  denote by
$\podd (n)$ the number of partitions of $n$ wherein odd
parts are distinct.
It is easy to see that
\begin{align*}
\sum_{n=0}^{\infty}\podd (n) q^n
=\prod_{n=1}^\infty \frac{1+q^{2n-1}}{1-q^{2n}}
=1+q+q^2+2q^3+3q^4+4q^5+5q^6+7q^7+\cdots.
\end{align*}

\begin{cor}\label{cor:4} For $n,k\geqslant 1$, there holds
\begin{align}
(-1)^{k-1}\sum_{j=0}^{k-1} (-1)^{j}\Big(\podd\big(n-j(2j+1)\big)-\podd\big(n-(j+1)(2j+1)\big)\Big)
\geqslant 0  \label{eq:poddneq}
\end{align}
with strict inequality if $n\geqslant (2k+1)k$. For example,
\begin{align*}
\podd (n)-\podd (n-1)-\podd (n-3)+\podd (n-6)\leqslant 0,  \\
\podd (n)-\podd (n-1)-\podd (n-3)+\podd (n-6)+\podd (n-10)-\podd (n-15)\geqslant 0.
\end{align*}
\end{cor}

A nice combinatorial proof of \eqref{eq:recur1} was given by Bressoud
and Zeilberger \cite{BZ}.  It would be interesting to find a
combinatorial proof of \eqref{eq:ineq1}, \eqref{eq:neq-overp} and \eqref{eq:poddneq}.
Moreover, Andrews and Merca \cite{AM}
found  a partition-theoretic interpretation of the truncated sum \eqref{eq:ineq1}.
It is still an open problem  to give partition  interpretations for our two truncated sums
in \eqref{eq:neq-overp} and \eqref{eq:poddneq}.
A combinatorial proof of \eqref{eq:neq-overp2} will be given in Section 4.

\section{Proof of Theorem \ref{thm:1}}
Generalizing Shanks' work \cite{Shanks,Shanks2},
Andrews~\cite[Lemma 2]{Andrews86} (see also Andrews, Goulden, and Jackson \cite[Theorem 1]{AGJ}) established the following identity
\begin{align}
\sum_{j=0}^n \frac{(b)_{j}(1-bq^{2j})(b/a)_j a^j q^{j^2}}{(1-b)(q)_j (aq)_j}
=\frac{(bq)_n}{(aq)_n}\sum_{j=0}^n \frac{(b/a)_j a^j q^{(n+1)j}}{(q)_j}. \label{eq:AGJ}
\end{align}
When $b=1$ and $a=-1$, the identity \eqref{eq:AGJ} reduces to

\begin{align}
1+2\sum_{j=1}^n (-1)^j q^{j^2}
=\sum_{j=0}^n (-1)^j \frac{(-1)_{j} (q)_n q^{(n+1)j}}{(q)_j (-q)_{n}}. \label{eq:agj-gauss}
\end{align}
By \eqref{eq:agj-gauss} and the $q$-binomial theorem (see \cite[Theorem 2.1]{Andrews98}), we have
\begin{align}
\frac{(-q)_\infty}{(q)_\infty}\bigg(1+2\sum_{j=1}^k (-1)^j q^{j^2}\bigg)
&=\sum_{j=0}^k (-1)^j \frac{(-1)_{j} (-q^{k+1})_\infty q^{(k+1)j}}{(q)_j (q^{k+1})_{\infty}} \nonumber\\
&=\sum_{j=0}^k (-1)^j \sum_{i=0}^{\infty}\frac{(-1)_{j} (-1)_i q^{(k+1)(i+j)}}{(q)_j (q)_i}.
\label{eq:overp}
\end{align}
After making a change of variable $i+j=n$ and reordering
the summation on the right-hand side of \eqref{eq:overp}, one  should then get  a double sum
$\sum_{n=0}^\infty \sum_{j=0}^{\min\{n,k\}}$. Since $1/(q)_m=0$ for $m=-1,-2,\dots$, one can write
this double sum as
\[
\sum_{n=0}^\infty \sum_{j=0}^k (-1)^j \frac{(-1)_{j} (-1)_{n-j} q^{(k+1)n}}{(q)_j (q)_{n-j}}.
\]
By induction on $k$, it is easy to see that, for $n\geqslant 1$,
\begin{align*}
\sum_{j=0}^k (-1)^j \frac{(-1)_{j} (-1)_{n-j}}{(q)_j (q)_{n-j}} =(-1)^k \frac{(-q)_{k} (-1)_{n-k}}{(1-q^n)(q)_{n-k-1} (q)_{k}}.
\end{align*}
Hence,  the right-hand side of \eqref{eq:overp} can be written as
\begin{align*}
&1+(-1)^k \sum_{n=1}^\infty \frac{(-q)_{k} (-1)_{n-k}q^{(k+1)n}}{(1-q^n)(q)_{n-k-1} (q)_{k}}\\
&=1+(-1)^k \sum_{n=k+1}^\infty \frac{(-q)_{k} (-1)_{n-k}q^{(k+1)n}}{(q)_n}{n-1\brack k},
\end{align*}
as desired.

\section{Proof of Corollary \ref{cor:2}}
By \eqref{eq:genovp} and \eqref{eq:ovp}, we see that the generating function for the sequence
$\{s_n\}_{n\geqslant 0}$, where
$$
s_n=(-1)^k \bigg(\overline{p}(n)+2\sum_{j=1}^k (-1)^{j}\overline{p}(n-j^2) \bigg),
$$
is given by
\begin{align}\label{eq:new}
(-1)^k+\sum_{n=k+1}^\infty \frac{(-q)_{k} (-1)_{n-k}q^{(k+1)n}}{(q)_n}{n-1\brack k}.
\end{align}
It follows that the  coefficient of $q^m$ in \eqref{eq:new} is  nonnegative for $m\geqslant 1$ and  positive
 for $m\geqslant (k+1)^2$. This completes the proof.

\section{A combinatorial proof of \eqref{eq:neq-overp2}}
Let $\mathscr P_n$ denote the set of all overpartitions of $n$. We now construct a mapping
$\phi\colon \mathscr P_n\to \mathscr P_{n-1}$ as follows: For any
$\lambda=(\lambda_1,\ldots,\lambda_k)\in \mathscr P_n$, let
\begin{align*}
\phi(\lambda)
=\begin{cases}(\lambda_1,\ldots,\lambda_{k-1}), &\text{if $\lambda_k=1$,}\\
(\lambda_1,\ldots,\lambda_{k-1},\lambda_k-1), &\text{if $\lambda_k\neq 1,\overline{1}$,} \\
(\lambda_1, \ldots, \lambda_{k-2}, \underbrace{\hat 1, \ldots, 1}_{\text{$\lambda_{k-1}$ $1$'s}} ),
&\text{if $\lambda_k=\overline{1}$,}
\end{cases}
\end{align*}
where $\hat 1=\overline 1$ if $\lambda_{k-1}$ is overlined and $\hat 1=1$ otherwise.

For example, for $n=4$,  the mapping $\phi$ gives
\begin{align*}
&4\mapsto 3,\ \overline{4}\mapsto\overline{3},\ (3,1)\mapsto 3,
\ (3,\overline{1})\mapsto (1,1,1),\ (\overline{3},\overline{1})\mapsto (\overline{1},1,1),
\ (\overline{3},1)\mapsto \overline{3},
\\
& (2, 2)\mapsto (2, 1),\ (\overline{2}, 2)\mapsto (\overline{2}, 1),
\ (2,1,1)\mapsto (2,1),
\ (2,\overline{1},1)\mapsto (2,\overline{1}),
\ (\overline{2},\overline{1},1)\mapsto (\overline{2},\overline{1}), \\
& (\overline{2},1,1)\mapsto (\overline{2},1),\ (1,1,1,1)\mapsto (1,1,1),
\ (\overline{1},1,1,1)\mapsto (\overline{1},1,1).
\end{align*}
It is easy to see that $1\leqslant |\phi^{-1}(\mu)|\leqslant 2$ for any $\mu\in \mathscr P_{n-1}$.
This proves that $\overline{p}(n)\leqslant 2\overline{p}(n-1)$.

\section{Proof of Theorem \ref{thm:3} and Corollary \ref{cor:4}}
In \cite{Shanks,Shanks2} Shanks proved that
\begin{align}
\sum_{j=0}^{n-1}q^{j(2j+1)}(1+q^{2j+1}) =\sum_{j=0}^{n-1}
\frac{(q;q^2)_j (q^2;q^2)_n q^{j(2n+1) }}{(q^2;q^2)_j (q;q^2)_n}.
\label{eq:shanks2}
\end{align}
By \eqref{eq:shanks2} (with $q$ replaced by $-q$) and
the $q$-binomial theorem (see \cite[Theorem 2.1]{Andrews98}), we have
\begin{align}
&\hskip -3mm \frac{(-q;q^2)_\infty}{(q^2;q^2)_\infty}
\sum_{j=0}^{k-1}(-1)^jq^{j(2j+1)}(1-q^{2j+1})\nonumber \\
&=\frac{(-q;q^2)_\infty}{(q^2;q^2)_\infty}
\sum_{j=0}^{k-1} (-1)^j \frac{(-q;q^2)_{j} (q^2;q^2)_k q^{(2k+1)j}}{(q^2;q^2)_j (-q;q^2)_{k}} \nonumber\\
&=\sum_{j=0}^{k-1} (-1)^j \frac{(-q;q^2)_{j} (-q^{2k+1};q^2)_\infty
q^{(2k+1)j}}{(q^2;q^2)_j (q^{2k+2};q^2)_{\infty}} \nonumber\\
&=\sum_{j=0}^{k-1} (-1)^j \sum_{i=0}^{\infty}\frac{(-q;q^2)_{j}
(-q^{-1};q^2)_i q^{(2k+1)(i+j)+i}}{(q^2;q^2)_j (q^2;q^2)_i}.
\label{eq:podd}
\end{align}
By induction on $k$, it is easy to see that, for $n\geqslant 1$,
\begin{align*}
\sum_{j=0}^{k-1}(-1)^j \frac{(-q;q^2)_{j}
(-q^{-1};q^2)_{n-j}q^{n-j}}{(q^2;q^2)_j (q^2;q^2)_{n-j}}
=(-1)^{k-1}\frac{(-q;q^2)_{k} (-q;q^2)_{n-k}
q^{n-k}}{(1-q^{2n})(q^2;q^2)_{n-k} (q^2;q^2)_{k-1}}.
\end{align*}
Hence, letting $i+j=n$, the right-hand side of \eqref{eq:podd} can be written as
\begin{align*}
&\hskip -3mm \sum_{n=0}^\infty
\sum_{j=0}^{k-1} (-1)^j \frac{(-q;q^2)_{j} (-q^{-1};q^2)_{n-j}q^{(2k+1)n+n-j}}{(q^2;q^2)_j (q^2;q^2)_{n-j}}\\
&=1+(-1)^{k-1} \sum_{n=1}^\infty \frac{(-q;q^2)_{k} (-q;q^2)_{n-k} q^{2(k+1)n-k}}
{(1-q^{2n})(q^2;q^2)_{n-k} (q^2;q^2)_{k-1}}.\\
&=1+(-1)^{k-1} \sum_{n=k}^\infty \frac{(-q;q^2)_{k} (-q;q^2)_{n-k}q^{2(k+1)n-k}}
{(q^2;q^2)_n}{n-1\brack k-1}_{q^2}.
\end{align*}
This proves Theorem \ref{thm:3}. The proof of Corollary \ref{cor:4} is similar to that of Corollary \ref{cor:2} and is
omitted here.

\section{Open problems}
In this section,
we propose a common   generalization  of \eqref{eq:ineq1}, \eqref{eq:neq-overp} and \eqref{eq:poddneq}.
Let $m,r$ be positive integers with $1\leqslant r\leqslant m/2$.
Consider
the generalized partition function $J_{m,r}(n)$ defined by
\begin{align}
\sum_{n=0}^\infty J_{m,r}(n) q^n=\frac{1}{(q^{r}, q^{m-r}, q^m;q^m)_\infty},\label{eq:jnr1}
\end{align}
where
$$
(a,b,c;q)_\infty=(a;q)_\infty(b;q)_\infty(c;q)_\infty.
$$
It is easy to see that
\begin{align*}
J_{2,1}(n)=\overline{p}(n),\quad
J_{3,1}(n)=p(n), \quad
J_{4,1}(n)=\podd (n).
\end{align*}
Moreover, if $r<m/2$, then $J_{m,r}(n)$ can be understood as the number of partitions of $n$ into parts
congruent to $0,\pm r$ modulo $m$. Now,
Jacobi's triple product identity implies (see \cite[p.~375]{HW}) that
\begin{align}
1+\sum_{j=1}^\infty (-1)^j (q^{j(mj+m-2r)/2}+q^{j(mj-m+2r)/2})=(q^{r}, q^{m-r}, q^m;q^m)_\infty \label{eq:jnr2}.
\end{align}
It follows from \eqref{eq:jnr1} and \eqref{eq:jnr2} that $J_{m,r}(n)$ satisfies the recurrence formula:
\begin{align*}
J_{m,r}(n)+\sum_{j=1}^{\infty} (-1)^{j}\Big(J_{m,r}\big(n-j(mj-m+2r)/2\big)
+J_{m,r}\big(n-j(mj+m-2r)/2\big)\Big)= 0,
\end{align*}
where $J_{m,r}(s)=0$ for all negative $s$.

\begin{conj}\label{conj:main}
For $m,n,k,r\geqslant 1$ with $r\leqslant  m/2$, there holds
\begin{align}
(-1)^{k-1}\sum_{j=0}^{k-1} (-1)^{j}\Big(J_{m,r}\big(n-j(mj+m-2r)/2\big)
-J_{m,r}\big(n-(j+1)(mj+2r)/2\big)\Big)
\geqslant 0 \label{eq:conjlast}
\end{align}
with strict inequality if $n\geqslant k(mk+m-2r)/2$.
\end{conj}

\begin{remark}
After we posted a preliminary version of this paper on arXiv, George E. Andrews  informed us that
our  Conjecture~\ref{conj:main} is effectively identical with Problem 2 in the final version of \cite{AM}.
\end{remark}

For $m=2$ and $r=1$,
the inequality \eqref{eq:conjlast} is equivalent to
\begin{align}
(-1)^{k-1} \bigg(\overline{p}(n)
+2\sum_{j=1}^{k-1} (-1)^{j}\overline{p}(n-j^2) \bigg)- \overline{p}(n-k^2)\geqslant 0  \label{eq:neq-overp-conj}
\end{align}
with strict inequality if $n\geqslant k^2$.
It is clear that \eqref{eq:neq-overp-conj} is stronger than the proved inequality
\eqref{eq:neq-overp} (with $k$ replaced by $k-1$).
By \eqref{eq:genovp} and \eqref{eq:ovp},  the generating function of the left-hand side of \eqref{eq:neq-overp-conj} is equal to
\begin{align}
&\hskip -2mm
(-1)^{k-1}\frac{(-q)_\infty}{(q)_\infty}\bigg(1+(-1)^{k} q^{k^2}
+2\sum_{j=1}^{k-1} (-1)^j q^{j^2}\bigg)\nonumber \\
&\qquad =(-1)^{k-1} +\sum_{n=k}^\infty
\frac{(-q)_{k-1} (-q)_{n-k}q^{kn}}{(q)_n}{n-1\brack k-1} \nonumber \\
&\hspace{2.5cm} -\sum_{n=k+1}^\infty\frac{(-q)_{k} (-q)_{n-k-1}q^{(k+1)n}}{(q)_n}{n-1\brack k}. \label{eq:newovp}
\end{align}
Therefore, the conjectured inequality \eqref{eq:neq-overp-conj} is equivalent to
\begin{conj} Let $k\geqslant 1$. Then
the coefficient of $q^m$ in the right-hand side of \eqref{eq:newovp} is nonnegative   for $m\geqslant 1$ and  positive for $m\geqslant k^2$.
\end{conj}

For $m=5$ and $r=1,2$,
the  inequality \eqref{eq:conjlast} may be written as
\begin{align}
(-1)^{k-1}\sum_{j=-k}^{k-1} (-1)^{j} J_{5,1}\bigl(n-j(5j+3)/2\bigr)
&\geqslant 0,   \label{eq:neq-RR1-conj} \\
(-1)^{k-1}\sum_{j=-k}^{k-1} (-1)^{j} J_{5,2}\bigl(n-j(5j+1)/2\bigr)
&\geqslant 0.  \label{eq:neq-RR2-conj}
\end{align}
By \eqref{eq:jnr1}, the conjectured inequalities
\eqref{eq:neq-RR1-conj} and \eqref{eq:neq-RR2-conj} are equivalent to
\begin{align}\label{eq:RR1}
&(-1)^k+\frac{(-1)^{k-1}}{(q,q^4,q^5;q^5)_\infty} \sum_{j=-k}^{k-1} (-1)^{j} q^{j(5j+3)/2}\in \N[[q]], \\
&(-1)^k+\frac{(-1)^{k-1}}{(q^2,q^3,q^5;q^5)_\infty} \sum_{j=-k}^{k-1} (-1)^{j} q^{j(5j+1)/2}\in \N[[q]].  \label{eq:RR2}
\end{align}
The two series $\sum_{j=-k}^{k-1}$  in \eqref{eq:RR1} and \eqref{eq:RR2} already appeared in
the works of Andrews \cite{Andrews83} and  Warnaar \cite{Warnaar}
(see also Chapman \cite{Chapman})
as  partial-sum analogues of Rogers-Ramanujan identities.  In particular, they obtained alternative expressions  of these series as double sums.
However, we have no idea  how to use  their formulas  to tackle the conjectures \eqref{eq:RR1} and \eqref{eq:RR2}.

Along the same line of thinking,  we consider the sequence
$\{t(n)\}_{n\geqslant 0}$ (see  A000716 in Sloane's database of integer sequences
\cite{Sloane}) defined by
\begin{align*}
\sum_{n=0}^{\infty} t(n) q^n
&=\frac{1}{(q)_\infty^3} \\[5pt] 
&=1+3q+9q^2+22q^3+51q^4+108q^5+221q^6+429q^7+810q^8+1479q^9+\cdots.
\end{align*}
Clearly, the number $t(n)$ counts partitions of $n$ into  $3$ kinds of parts. Now, invoking
the identity of Jacobi~\cite[p.~377]{HW}:
\begin{align*}
{(q)_\infty^3}=\sum_{j=0}^\infty (-1)^j (2j+1)q^{j(j+1)/2},
\end{align*}
we derive  the recurrence formula:
\begin{align*}
\sum_{j=0}^\infty (-1)^j(2j+1)t(n-j(j+1)/2)=0,
\end{align*}
where $t(m)=0$ for all negative $m$.

We end the paper with the following conjecture:
\begin{conj}For $n,k\geqslant 1$, there holds
\begin{align*}
(-1)^k\sum_{j=0}^k (-1)^j(2j+1)t(n-j(j+1)/2) \geqslant 0
\end{align*}
with strict inequality if $n\geqslant (k+1)(k+2)/2$. For example,
\begin{align*}
t(n)-3t(n-1)\leqslant 0,  \\
t(n)-3t(n-1)+5t(n-3)\geqslant 0, \\
t(n)-3t(n-1)+5t(n-3)-7t(n-6)\leqslant 0.
\end{align*}
\end{conj}

\vskip 5mm \noindent{\bf Acknowledgments.} We are grateful to George E. Andrews for
helpful comments on a previous version of this paper and for sending us the final version of \cite{AM}.
We also thank the two referees for valuable suggestions.
This work was  partially supported by CMIRA COOPERA 2012 de  la R\'egion Rh\^one-Alpes.

\end{document}